\documentclass[letterpaper,10pt,draftcls,onecolumn]{IEEEtran}
\usepackage{amsmath,amssymb, amsfonts}
\usepackage{hyperref}
\usepackage{graphicx}
\usepackage{wrapfig,lipsum}
\usepackage{mathrsfs}
\usepackage{xcolor}
\usepackage{tikz}

\newtheorem{thm}{Theorem}[section]

\newtheorem{lem}[thm]{Lemma}

\newtheorem{probl}{Problem}

\newtheorem{rem}[thm]{Remark}
\newtheorem{defn}[thm]{Definition}

\newtheorem{remark}{Remark}

\newcommand{\R}{\mathbb{R}}

\newcommand{\N}{\mathbb{N}}

\newcommand{\bit}{\begin{itemize}}
\newcommand{\eit}{\end{itemize}}
\newcommand{\ben}{\begin{enumerate}}
\newcommand{\een}{\end{enumerate}}

\usepackage{mathrsfs}
\usepackage{authblk}

\newcommand{\lp}{\left (}
\newcommand{\rp}{\right )}

\DeclareMathOperator{\tr}{tr}
\newcommand{\Ex}[1]{\mathbb{E}\left[ #1\right]}
\newcommand{\Var}[1]{\rm{Var}\left[ #1\right]}
\newcommand{\Prob}[1]{\mathbb{P}\left[ #1\right]}

\newcommand{\nodeset}{V}
\newcommand{\nodenum}{n}
\newcommand{\edgeset}{E}
\newcommand{\graph}{G}
\newcommand{\randset}{\mathcal{G}}
\newcommand{\lapset}{\mathcal{L}}
\newcommand{\union}{U}

\newcommand{\unionset}{\mathcal{U}}

\usepackage{mathtools}
\mathtoolsset{showonlyrefs=true}

\newcommand{\blue}[1]{#1}

\newtheorem{proposition}{Proposition}
\newtheorem{example}{Example}

\makeatletter            
\let\@fnsymbol\@arabic 
\makeatother             

\title{\blue{Non-Asymptotic} Connectivity of Random Graphs and Their Unions}
\author{Beth Bjorkman$^{*}$\thanks{$^{*}$Iowa State University Department of Mathematics, \texttt{bjorkman@iastate.edu}} \and 
Matthew Hale$^{\dagger}$\thanks{$^{\dagger}$Department of Mechanical and Aerospace Engineering, University of Florida, \texttt{matthewhale@ufl.edu}} \and  
Thomas Lamkin$^{\ddagger}$\thanks{$^{\ddagger}$lamkintd@miamioh.edu} \and 
Benjamin Robinson$^{**}$\thanks{$^{**}$Air Force Research Laboratory Sensors Directorate} \and Craig Thompson$^{\dagger\dagger}$\thanks{$^{\dagger\dagger}$University of Arizona Program in Applied Mathematics, \texttt{craigthompson@math.arizona.edu}}}

\begin{document}
\maketitle

\begin{abstract}
Graph-theoretic methods have seen wide use throughout the
literature on
multi-agent control and optimization. When
communications are intermittent and unpredictable, such
networks have been modeled using random communication
graphs. When graphs are time-varying, it is common
to assume that their unions are connected over time, yet,
to the best of our knowledge, there are not results that
determine the number of finite-size random graphs needed to attain
a connected union. Therefore, this paper bounds the probability
that individual random graphs are connected and bounds
the same probability for connectedness of unions of random graphs. 
The random graph model used is a generalization of the classic
Erd\H{o}s-R\'{e}nyi model which allows some edges
never to appear. Numerical results are presented to illustrate
the analytical developments made. 
\end{abstract}

\section{Introduction}\label{sec:intro}
Multi-agent systems have been studied in a number of applications, including sensor networks~\cite{cortes02}, 
communications~\cite{kelly98}, and smart power grids~\cite{caron10}. 
Across these applications, the agents in a network and their associated communications are often abstractly represented as graphs \cite{mesbahi10}. In general, graph-theoretic methods in multi-agent systems represent each agent as a node in a graph and each communication link as an edge, and multi-agent coordination algorithms have been developed for both static and time-varying graphs \cite[Chapter 1.4]{mesbahi10}. 

Time-varying random graphs in particular have been used to model communications which are unreliable and intermittent due to interference and poor channel quality \cite[Chapter 5]{mesbahi10}. Such graphs have seen use in a number of multi-agent settings. For example, distributed agreement problems over random graphs are studied in \cite{touri11} and \cite{hatano05}, while optimization over random graphs was explored in \cite{lobel11}. The work in \cite{yazicioglu15} provides a means to modify random graphs to make them robust to network failures, and \cite{lewis13} discusses general properties of random graphs as they pertain to multi-agent systems.  A broad survey of graph-theoretic results for control can be found in \cite{mesbahi10}, and well-known graph-theoretic results in optimization include \cite{blondel05,nedic09,zhu12}. 

When time-varying graphs (random or not) are used, a common assumption is that the unions of these graphs are connected over intervals of some finite length, i.e., the graph containing all edges present over time is itself a connected graph. A partial sampling of works using this assumption (or a related variant) includes \cite{moreau05,ren05,tseng90,tseng91,ren05b,chen12,kia15,blondel05,jadbabaie03,touri09,nedic10,nedic09,olshevsky11,zhu12,feyzmahdavian14,nedic07}.  In addition, some works derive convergence rates or other results that explicitly use the length of such intervals, including \cite{feyzmahdavian14,tseng90,nedic07,tseng91,chen12,blondel05,touri09,olshevsky11,nedic09}. In applying these results, one may wish to determine the time needed for the system to attain a connected union graph. To the best of our knowledge, no study has been undertaken that addresses this problem for unions of finite-size random graphs, despite their frequent use in multi-agent systems. 

 Owing to the success of Erd\H{o}s-R\'{e}nyi graphs in modeling some time-varying multi-agent communications \cite[Chapter 5]{mesbahi10}, we consider unions of random graphs generated by a generalization of the Erd\H{o}s-R\'{e}nyi model.
In particular, we examine the connectedness of such
graphs and their unions. Formally, this paper 
first lower-bounds the number of graphs needed to make
their ``expected union'' connected (in a precise sense). 
It then finds a lower bound for the probability that a union of random graphs is connected as a function
of graph model parameters.

Our results use spectral properties of the Laplacian of a union of random graphs. We utilize the first order statistic of the set of non-trivial eigenvalues to bound the expected value of the Laplacian's second-smallest eigenvalue, called the \emph{algebraic connectivity} \cite{fiedler73} of the underlying union graph.  This bound in turn enables  a lower bound on the probability of the algebraic connectivity exceeding some given threshold. 

Random graphs' Laplacians are random matrices, and thus our
approach relies on the spectral properties of random matrices
of a particular form.  One common approach to analyzing the spectra of random matrices is to let the dimension of the matrix grow arbitrarily large \cite{diaconis94,furedi81,wigner58}; the work in \cite{coja07} considers similar asymptotic results focused specifically on Laplacians of random graphs. 
For random graphs, a common approach is to derive results in which the size of the graph grows arbitrarily large, and doing so enables results that hold for \emph{almost all} graphs \cite{bollobas01}. 
\blue{
In particular, seminal work by Erd\H{o}s and R\'{e}nyi showed that
an Erd\H{o}s-R\'{e}nyi graph on~$n$ nodes
with edge probability~$p$ is almost surely connected as~$n \to \infty$
if~$p \geq c\frac{\log n}{n}$ and~$c > 1$. 
}
While there is clear theoretical appeal to such results, our focus on multi-agent systems leads us to consider non-asymptotic results precisely because such systems are very often comprised by a fixed, finite number of agents. 
\blue{We therefore study and characterize random graphs with a fixed, finite number of nodes, i.e., with~$n$ fixed, and we will  develop a novel statistical approach to characterize such graphs.} 

In addition, while \blue{references~\cite{erdos59} and~\cite{erdos60}, and many other} 
works on random graphs consider edge probabilities that bear some known relationship to the number of nodes in a graph \cite{krivelevich03,tran13}, we do not do so here. 
The use of random graphs to model multi-agent communications is inspired by applications in which poor channel quality, interference, and other factors make communications unreliable. In such cases, the probability of a communication link being active may not bear any known relationship to the size of the network. Accordingly, we proceed with edge probabilities and network sizes that are fixed and not assumed to be related. \blue{In particular, we require only~$p \in (0, 1)$ and avoid
any assumption that~${p \sim \frac{\log n}{n}}$ or that~$p$
is any known function of~$n$, in contrast to existing results such as~\cite{erdos60} and~\cite{bollobas01}}. 

\blue{
Related work in~\cite{benjamini08} explores the properties
of isoperimetric constants of random graphs, which
are known to be related to algebraic connectivity
via Cheeger's Inequality~\cite{cheeger69}. However, the random
graph model in~\cite{benjamini08} always adds one new edge at every
time step, while the unions of random graphs we consider
can add zero, one, or many new edges at every time step.
Accordingly, the conclusions reached in~\cite{benjamini08} do not apply
to the random graph model we consider. 
Developments in~\cite{cooper07} assess connectivity of 
the Erd\H{o}s-R\'{e}nyi random graph model
through analyzing their cover times. However,
the results in~\cite{cooper07} are asymptotic, in that
they consider~$n \to \infty$, and thus their
conclusions do not readily apply to the finite-sized
graphs we study. 
}

\blue{
Additional work in~\cite[Section V]{pirani15}
studies eigenvalues of Laplacians of Erd\H{o}s-R\'{e}nyi graphs
as we do in this work. 
In particular, results in that section study
the algebraic connectivity of Erd\H{o}s-R\'{e}nyi graphs,
though they do so in the asymptotic regime (i.e.,
as~$n \to \infty$). Our results differ from that work
because they apply for arbitrary, fixed values of~$n$,
rather than to the case in which~$n \to \infty$. 
}

The rest of the paper is organized as follows.  Section~\ref{sec:review} reviews the required elements of graph theory and gives a formal statement for the problem that is the focus of this paper. Section~\ref{sec:firstorder} computes the first order statistic of the eigenvalues of random graph Laplacians and certain statistical properties of the eigenvalues to enable the results of Section~\ref{sec:algconn}. Section~\ref{sec:algconn} then presents the main results of the paper and solves the problem stated in Section \ref{sec:review}.
Section~\ref{sec:specific} next provides numerical solutions to several instantiations of the problem studied. Finally, Section~\ref{sec:conclusion} provides concluding remarks.

\section{Review of Graph Theory and Problem Statement} \label{sec:review}
In this section, we review the required elements of graph theory. We begin with basic definitions, including the definition of algebraic connectivity, and then review the Erd\H{o}s-R\'{e}nyi model for random graphs and introduce a generalization that we use. Throughout this paper, all uses of the phrase ``random graphs'' refer to Erd\H{o}s-R\'{e}nyi-type graphs. Then we formally state the statistical graph
connectivity problem solved in this paper. Below, we use the notation $[n] := \{1, \ldots, n\}$ for any $n \in \N$ and $|\cdot|$ to represent the cardinality of a set.  

\subsection{Basic Graph Theory}
A \emph{graph} $\graph=(\nodeset,\edgeset)$ is defined over a set of nodes, denoted $\nodeset$, and describes connections between these nodes in the form of edges, which are contained in an edge set $\edgeset$. For $\nodenum$ nodes, $\nodenum \in \N$, the elements of $\nodeset$ are indexed over $[\nodenum]$. The set of edges in the graph is a subset $\edgeset \subseteq \nodeset \times \nodeset$, where a pair $(i, j) \in \edgeset$ if nodes $i$ and $j$ share a connection, and $(i, j) \not\in \edgeset$ if they do not. This paper considers graphs which are \emph{undirected}, meaning an edge $(i, j)$ is not distinguished from an edge $(j, i)$, and \emph{simple}, so that $(i, i) \not\in \edgeset$ for all $i$. The \emph{degree} of node $i\in\nodeset$ is defined as 
\begin{equation}
d_i=\big|\{ j\in\nodeset \mid (i,j)\in\edgeset \}\big|, 
\end{equation}
which is equal to 
the number of vertices sharing an edge with
vertex~$i$. One main focus of this paper is \emph{connected} graphs.
\begin{defn}[E.g., \cite{godsil01}]
A graph $\graph$ is called \emph{connected} if for all $i \in [\nodenum]$ and $j \in [\nodenum]$, $i \neq j$, there is a sequence of edges one can traverse from node $i$ to node $j$, i.e., there is a sequence of vertices $\{i_{\ell}\}_{\ell=1}^{k}$ 
such that $\edgeset$ contains all of the edges
\begin{equation}
(i, {i_1}), ({i_1}, {i_2}), ({i_2}, {i_3}), \ldots, ({i_{k-1}}, {i_k}), ({i_k}, j).
\end{equation}
\hfill $\triangle$
\end{defn}

The results of this paper are developed in terms of graph Laplacians, which are defined in terms of the adjacency and degree matrices of a graph. The \emph{adjacency matrix} ${A(\graph) \in \R^{\nodenum \times \nodenum}}$ associated with the graph $\graph$ is defined element-wise as 
\begin{equation}
(A(\graph))_{ij} = \begin{cases} 1 & (i, j) \in \edgeset \\
                                 0 & \textnormal{otherwise}
                                 \end{cases}.
\end{equation}
 When there is no ambiguity, we will simply denote $A(\graph)$ by $A$. Because we consider undirected graphs, $A$ is symmetric by definition. The degree matrix $D(\graph) \in \R^{\nodenum \times \nodenum}$ associated with a graph $\graph$ is the diagonal matrix  
 ${D(\graph) = \textnormal{diag}(d_1, d_2, \ldots, d_n)}$, which we will denote $D$ when $\graph$ is clear from context. By definition, $D$ is also symmetric.  The \emph{Laplacian} of a graph $\graph$ is then defined as ${L(\graph) = D(\graph) - A(\graph)}$, which will be written simply as $L$ when $\graph$ is unambiguous. 
 
 The results of this paper rely in particular on spectral properties of $L$. Letting $\lambda_k(\cdot)$ denote the $k^{\text{th}}$ smallest eigenvalue of a matrix, it is known that $\lambda_1(L) = 0$ for all graph Laplacians \cite{mesbahi10}, and thus we have
\begin{equation} \label{eq:lambdai}
0 = \lambda_1(L) \leq \lambda_2(L) \leq \cdots \leq \lambda_{\nodenum}(L).
\end{equation}
The value of $\lambda_2(L)$ is central to the work in this paper and some other works in graph theory, and it gives rise to the following definition. 

\begin{defn}[\cite{fiedler73}] \label{def:lambda2}
The \emph{algebraic connectivity} of a graph $\graph$ is the second smallest eigenvalue of its Laplacian, $\lambda_2(L)$, and $\graph$ is connected if and only if $\lambda_2(L) > 0$. \hfill $\triangle$
\end{defn}

This paper is dedicated to studying the statistical properties of $\lambda_2=\lambda_2(L)$ for 
random graphs and unions of random graphs. Toward doing so, we now review the basics of the random graph model we study. 

\subsection{Random Graphs}
Several well-known random graph models exist in the literature \cite{erdos59,watts98}, and Erd\H{o}s-R\'{e}nyi graphs in particular have been successfully used in the multi-agent systems literature. Erd\H{o}s-R\'{e}nyi graphs can model, for example, unreliable, intermittent and time-varying communications in multi-agent networks \cite{mesbahi10}, and we therefore consider the Erd\H{o}s-R\'{e}nyi model in this paper. Under this model, a graph on $\nodenum$ vertices contains each admissible edge with some fixed \emph{edge probability} $p \in (0, 1)$. Therefore, for each $i \in [\nodenum]$ and $j \in [\nodenum]$ with $i \neq j$, an Erd\H{o}s-R\'{e}nyi graph satisfies 
\begin{equation}
\Prob{(i, j) \in \edgeset} = p \,\, \textnormal{ and } \,\, \Prob{(i, j) \not\in E} = 1 - p.
\end{equation}
We denote the sample space of all Erd\H{o}s-R\'{e}nyi graphs on $\nodenum$ nodes with edge probability $p$ by $\randset(n, p)$, and we denote the set of Laplacians of all such graphs by $\lapset(n, p)$. 
We will primarily study the following variant of the Erd\H{o}s-R\'{e}nyi
model, which contains the typical
Erd\H{o}s-R\'{e}nyi formulation as a special case.  
\begin{defn}
Let $G = \big(V(G), E(G)\big)$ be a connected graph  on $n$ vertices. An \emph{Erd\H{o}s-R\'{e}nyi graph corresponding to $G$} is a random graph on $n$ vertices with an edge set $E$ which satisfies
\begin{eqnarray*}
\Prob{(i,j)\in \edgeset}= \left\{\begin{array}{rl}
p & (i,j)\in \edgeset(G)\\
0 & (i,j)\not\in\edgeset(G)
 \end{array} \right.. \end{eqnarray*}
Correspondingly,
\begin{eqnarray*}
\Prob{(i,j)\not\in \edgeset}= \left\{\begin{array}{rl}
1-p & (i,j)\in \edgeset(G)\\
1 & (i,j)\not\in\edgeset(G)
 \end{array} \right.. \end{eqnarray*}
 We denote the sample space of all Erd\H{o}s-R\'{e}nyi graphs corresponding to $G$ by $\randset(G,p)$.
We note that $\mathcal{G}(n,p)= \mathcal{G}(K_n,p)$, where $K_n$ represents the complete graph on $n$ vertices.
\hfill $\triangle$
\end{defn}


Some well-known results in the graph theory literature assume that $p$ has some known relationship to $\nodenum$ \cite{bollobas98}, or else that the number of edges in a random graph has some relationship to the number of nodes in the graph \cite{erdos60}. While the theoretical utility of these relationships is certainly clear from those works, these relationships will often not hold in multi-agent systems simply because their communications are affected by a wide variety of external factors. We therefore proceed with a value of $p \in (0, 1)$ that is not assumed to have any relationship to the value of $\nodenum$.  In the study of multi-agent systems, it is also common for algorithms and results to be stated in terms of unions of graphs, which we define now.  
\begin{defn}
For a collection of graphs ${\{\graph_k = (V, E_k)\}_{k=1}^{T}}$ defined on the same node set $\nodeset$, the union of these graphs, denoted $\union_{T}$, is defined as 
${\union_{T} := \cup_{k=1}^{T} G_k = \left(V, \cup_{k=1}^{T} E_k\right)}$,
i.e., the \emph{union graph} $\union_T$ contains all edges from all $T$ graphs in the union. \hfill $\triangle$
\end{defn}

\subsection{Problem Statement}
A common requirement in some multi-agent systems is that the communication graphs in a network form a connected union graph over intervals of some fixed length. To help determine when this occurs with random interactions, 
we solve the following problem in this paper.

\begin{probl} \label{prob:prob}
Given a connected graph $G$ \blue{on a fixed, finite number of nodes~$n \in \N$}
and 
\blue{a fixed, arbitrary} probability $p \in (0, 1)$, 
find a lower bound for the probability that a given graph 
$\hat{G} \in \mathcal{G}(G,p)$ is connected. That is, find a lower bound for $\mathbb{P}[\lambda_2(\hat{G}) > 0]$. \hfill $\diamond$
\end{probl}

%
 
While this problem is defined in terms of a single random graph, we can use it to derive results for a union of random graphs because such unions are themselves equivalent to single random graphs, albeit with a different edge probability.

\begin{lem} \label{lem:unioneq}
Let $\unionset_T(G,p)$ denote the set of all unions of $T$ random graphs on $\nodenum$ 
nodes with edge probability $p$, i.e., $ \unionset_T(G, p) := \left\{\cup_{k=1}^{T} \graph_k \mid \graph_k \in \randset(G, p)\right\}. $
Then
\begin{eqnarray*}
\unionset_T(G, p) = \randset\lp G, 1 - (1 - p)^T \rp.
\end{eqnarray*}
\end{lem}
\emph{Proof:} 
Fix any edge $(i,j)\in E(G)$. Then the edge $(i,j)$ is absent  in $\unionset_T(G,p)$ only if it is absent in all $T$ graphs that comprise $\unionset_T (G,p)$. Observe that the edge between $i$ and $j$ is absent in each $G_k$ with probability $(1-p)^T$. Then edge $(i,j)$ is present in $\unionset_T (G,p)$ with probability $1 - (1 - p)^T$.
\hfill $\blacksquare$

With Lemma~\ref{lem:unioneq}, results pertaining to individual random graphs can easily be
applied to unions of such graphs.

\subsection*{\blue{D. Differences with Asymptotic Analyses}}
\blue{
In this subsection, we highlight the differences between
what is asked in Problem~$1$ and classical results
in the study of random graphs. In particular, the asymptotic connectivity 
of Erd\H{o}s-R\'{e}nyi graphs is well-known. To highlight the key
differences with that result, we reproduce it here.
\begin{proposition}[From~\cite{erdos60}] \label{prop:asymptotic}
A graph~$G \in \mathcal{G}(n, p)$ is almost surely
connected if~$p \geq c\frac{\log n}{n}$
and~$c > 1$. Explicitly, 
\begin{equation}
\lim_{n \to \infty} \mathbb{P}\Big[G \in \mathcal{G}(n, p) \textnormal{ is connected}\Big] = 1
\end{equation}
if~$p = c\frac{\log n}{n}$ and~$c > 1$. \hfill $\blacklozenge$
\end{proposition}
}

\blue{
The notion of ``almost surely'' used in the statement
of Proposition~\ref{prop:asymptotic} dates back
to seminal work of Erd\H{o}s and R\'{e}nyi
and refers to properties that hold ``with probability
tending to~$1$ for~$n \to +\infty$''~\cite[Page 3]{erdos60}. 
That is, drawing conclusions about random graphs that hold ``almost surely'' 
requires taking the limit as~$n \to \infty$. 
This notion is sometimes
called ``asymptotically almost surely'', e.g.,
in~\cite{pirani15},
to reflect the inherently asymptotic nature of its conclusions. 
}

\blue{
The focus on asymptotics 
is the first essential difference between
what Problem~1 seeks to solve and existing
results: Problem~1 considers networks
of a fixed but arbitrary number of nodes~$n \in \N$, while
Proposition~\ref{prop:asymptotic} assesses connectivity
as~${n \to \infty}$. The finiteness of~$n$ in Problem~1
precludes the use of Proposition~\ref{prop:asymptotic} to solve it. 
The second key difference is that Problem~1
does not permit the assumption that~$p$ is a specified function
of~$n$. 
Instead, in Problem~1,~$p \in (0, 1)$
is fixed but arbitrary. 
Thus, any connectivity conclusions drawn
from Proposition~\ref{prop:asymptotic} are necessarily
inapplicable to Problem~1. 
}

\blue{
We illustrate the need for new techniques to solve
Problem~1 with an example below.
In it, we show that the setting of Problem~1
requires us to consider graphs that are neither
connected with probability~$1$, nor disconnected
with probability~$1$. 
}

\begin{example}\label{ex:graphs}
\textcolor{blue}{Fix~${n = 3}$ nodes and
consider the underlying
graph~${G = K_3}$, i.e., the complete graph on~$3$ nodes. 
Set~${p = \frac{1}{2} = 1.365\frac{\log 3}{3}}$.
If we take~${n \to \infty}$, then this choice of~$p$
satisfies the condition for connectivity in
Proposition~\ref{prop:asymptotic}, but the probability
that a graph in~$\mathcal{G}(3, \frac{1}{2})$ is connected
is not near~$1$, as follows.
We consider graphs in~${\mathcal{G}(K_3, \frac{1}{2})
= \mathcal{G}(3, \frac{1}{2})}$, which are 
classical Erd\H{o}s-R\'{e}nyi graphs. 
With~$n=3$ fixed, 
there are eight possible
graphs in~$\mathcal{G}(3, \frac{1}{2})$, shown in 
Figure~\ref{fig:eight}.}

\begin{figure}
\centering
\includegraphics[scale=0.8]{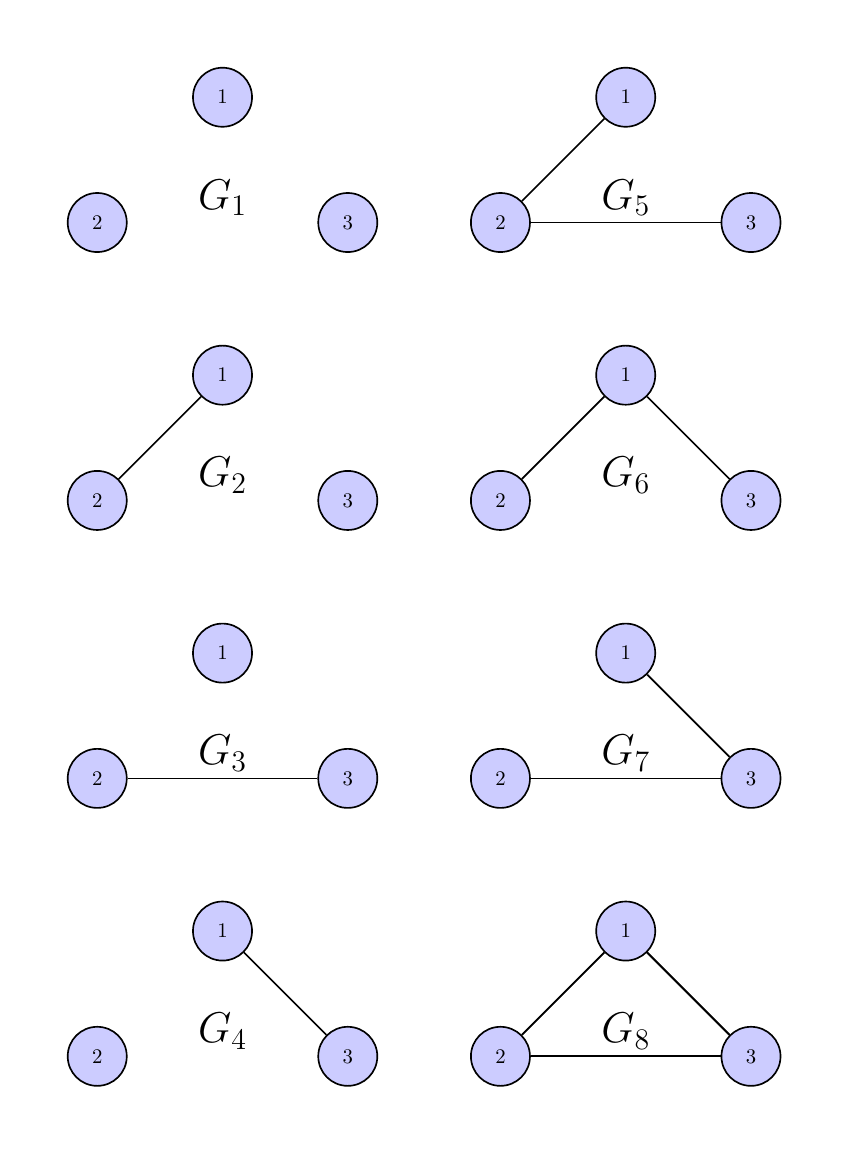}
\caption{\blue{The set of all eight graphs in~$\mathcal{G}(K_3, \frac{1}{2})= \mathcal{G}(3, \frac{1}{2})$. All graphs are equiprobable
because~${p = \frac{1}{2}}$, and thus half of all
graphs in~$\mathcal{G}(3, \frac{1}{2})$ are connected.}}
\label{fig:eight}
\end{figure}

\blue{
Graphs~$G_1$, $G_2$, $G_3$, and~$G_4$ are not connected,
while graphs~$G_5$, $G_6$, $G_7$, and~$G_8$ are connected.
Because~$p = \frac{1}{2}$, all eight graphs are
equiprobable and thus
\begin{equation}
\mathbb{P}\Big[G \in \mathcal{G}\big(3, \frac{1}{2}\big) \textnormal{ is connected}\Big] = \frac{1}{2}. 
\end{equation}
Then a random graph~$G \in \mathcal{G}(3, \frac{1}{2})$
is neither connected nor disconnected
with probability~$1$.
} 
\hfill $\lozenge$
\end{example}

\blue{
In Example~1, the lack of either connectedness or 
disconnectedness with probability~$1$ is not
captured by classical results such as Proposition~\ref{prop:asymptotic}, and
new analyses are required to account for and
quantify this behavior. That is the subject
of this paper. 
}

Section~\ref{sec:firstorder} next provides theoretical developments that will enable the solution to 
Problem~\ref{prob:prob} in Section~\ref{sec:algconn}.

\section{Order Statistics over Random Eigenvalues} \label{sec:firstorder}
Solving Problem~\ref{prob:prob} will require
bounds on the expected algebraic connectivity of a random graph
of interest. This section derives such bounds in terms of known
quantities by relating the expected algebraic connectivity
of a random graph to other spectral properties.

Let~$G$ be a connected graph with~$m$ edges and~$n$ vertices labeled by~$[n]$. Let~$d_i$ be the degree of vertex~$i$ in~$G$.
Let~$\hat{G} \in \mathcal{G}(G, p)$ be a random graph
with Laplacian~$L$ whose eigenvalues are denoted~$\lambda_1(\hat{G})$
through~$\lambda_n(\hat{G})$.
Let $\ell$ be a random
variable equal to~$\lambda_i$ with probability $\frac{1}{n-1}$ for $i\geq 2$. Take $N$ independent samples of $\ell$, denoted $\ell_1,\ldots,\ell_N$, and place the results in ascending order to create the list
\begin{eqnarray*}\ell_{1:N},\ell_{2:N},\ldots,\ell_{N:N},\end{eqnarray*}
where~$\ell_{j:N}$ denotes the~$j^{th}$ largest sample
out of the~$N$ total samples. 
Then $\ell_{1:N}$ is by definition the first order statistic for
the samples~$\{\ell_{j}\}_{j=1}^{N}$.
Analyzing these order statistics for $\hat{G}$,
we obtain the following relationship
to the expected value of~$\lambda_2(\hat{G})$. 

\begin{lem}\label{lem:l1Nupperbound} Consider a connected graph $G$ with $m$ edges and vertex set indexed over $[n]$.
For~$\hat{G} \in \mathcal{G}(G, p)$ we have
\begin{equation}
\Ex{\ell_{1:N}} \leq\Ex{\lambda_2(\hat{G})}\lp\! 1-\lp \frac{n-2}{n-1}\rp^{N-1}\rp+ \frac{2mp}{n-2}\lp \frac{n-2}{n-1}\rp^N.
\end{equation}		
\end{lem}
\emph{Proof:} 
Let $r=\frac{1}{n-1}$. Then
\begin{align}
\Ex{\ell_{1:N}} &= \Ex{\ell_{1:N} \mid \exists k : \ell_{k} = \lambda_2}\Prob{\exists k : \ell_k = \lambda_2} + \Ex{\ell_{1:N} \mid \not\exists k : \ell_{k} = \lambda_2}\Prob{\not\exists k : \ell_k = \lambda_2} \\
			&= \Ex{\ell_{1:N} \mid \exists k : \ell_{k} = \lambda_2}\lp1-\lp1-r\rp^N\rp + \Ex{\ell_{1:N} \mid \not\exists k : \ell_{k} \neq \lambda_2}\lp 1-r\rp^N.\label{eq:lem31runback}
\end{align}
Observe that
\begin{equation} \label{eq:exp1}
\Ex{\ell_{1:N} \mid \exists k : \ell_{k} = \lambda_2} = \mathbb{E}[\lambda_2].
\end{equation}
 
Furthermore, by definition we have
\begin{align}
\Ex{\ell_{1:N} \mid \not\exists k : \ell_{k} \neq \lambda_2} &:= 
\Ex{\min_{j \in [N]} \ell_j \mid \not\exists k : \ell_{k} \neq \lambda_2} \\
&\leq \Ex{\ell_j, \,\, j \textnormal{ arbitrary} \mid \not\exists k : \ell_{k} \neq \lambda_2} \\
&= \frac{1}{n-2}\sum_{i=3}^{n} \Ex{\lambda_i} \\
&= \frac{1}{n-2}(\Ex{\tr L} - \Ex{\lambda_2}),\label{eq:exp2}
\end{align}
where the last line follows from the fact that~$\Ex{\lambda_1} = 0$
for all graphs~$G$. 

Substituting Equations~\eqref{eq:exp1} and~\eqref{eq:exp2}
into Equation~\eqref{eq:lem31runback} gives
\begin{equation}
\Ex{\ell_{1:N}} \leq \Ex{\lambda_2}\lp1-\lp1-r\rp^N\rp + \lp\frac{\Ex{\tr L}-\Ex{\lambda_2}}{n-2}\rp\lp 1-r\rp^N.
\end{equation}
Next, we note that~${\Ex{\tr L} = \sum_{i=1}^{n} pd_i}$ by construction
of the Laplacian.
Then we have
\begin{equation}
\Ex{\ell_{1:N}} \leq \Ex{\lambda_2}\lp1-\lp1-r\rp^N-\frac{1}{n-2}\lp 1-r\rp^N\rp + \frac{\sum_{i=1}^n pd_i}{n-2}\lp 1-r \rp^N 
\end{equation}
and the lemma follows by factoring and the Handshake Lemma ($\sum_{i=1}^n d_i = 2m$).
\hfill $\blacksquare$

Having put~$\Ex{\lambda_2(\hat{G})}$ in terms
of~$\Ex{\ell_{1:N}}$, we next focus on bounding~$\Ex{\ell_{1:N}}$ 
in terms of known quantities. Toward doing so, we state
the following lemma. 

\begin{lem}[Equation (4), \cite{arnold79}]\label{lem:hale6}
Let $X_1, X_2, \ldots, X_m$ be jointly distributed with common mean $\mu$ and variance $\sigma^2$. Then the $k^{th}$ order statistic of this collection, denoted $X_{k:m}$, has expectation bounded according to 
\begin{eqnarray*}
\mu-\sigma \sqrt{\frac{m-k}{k}}\leq \Ex{X_{k:m}}.
\end{eqnarray*}
\hfill $\blacksquare$
\end{lem}

\blue{
The bound in Lemma~3.2 was shown in~\cite{arnold79} to be tight
for collections of random variables with identical
means and variances. 
}
\blue{Because each~$\ell_k$ has the same
mean and variance, Lemma~3.2 is an attractive
choice in computing the required bound on~$\mathbb{E}[\ell_{1:N}]$.
Therefore,  } we now work to apply Lemma~3.2 to $\ell_{1:N}$. 

\begin{lem}\label{lem:musigma}
Consider a connected graph~$G = (V, E)$ 
on~$n$ vertices with~$m$ edges so that 
vertex~$i$ has degree~$d_i$, and suppose that $\{\ell_k\}_{k \in [N]}$ are as before.
Then,
for~$\hat{G} \in \mathcal{G}(G, p)$,
we find that each~$\ell_k$ has mean~$\mu$ and
variance~$\sigma^2$ given by
\begin{equation}
\mu=\frac{2mp}{n-1}
\end{equation}
and
\begin{equation}
\sigma^2 \!=\! \frac{p}{(n-1)^2}\lp \! \lp n-1\rp\lp 2m(2-p)+p\sum_{i=1}^n  d_i^2\rp -4m^2p\rp.
\end{equation}
\end{lem}
\emph{Proof:} 
Observe that
\begin{equation}
\mu = \Ex{\frac{1}{n-1}\tr(L)} = \frac{2mp}{n-1}.
\end{equation}
Next, 
set $r=\frac{1}{n-1}$. 
Let~$L$ be the Laplacian of~$\hat{G}$, and
denote its~$i^{th}j^{th}$ off-diagonal entry
by~$X_{ij}$. If~$(i, j) \not\in E$, then~$X_{ij} = 0$.
If~$(i, j) \in E$, then
\begin{equation}
X_{ij} = \begin{cases} 1 & \textnormal{with probability~$p$} \\
 0 & \textnormal{with probability~$1-p$} \end{cases},
\end{equation}
i.e.,~$X_{ij}$ is a Bernoulli random variable
when~$(i, j) \in E$.  
Then $X_{ii}=\sum_{1\leq j \leq n, j\neq i} X_{ij}$.  
By definition we have
\begin{eqnarray}
\sigma^2	&=&\Ex{\ell_{k}^2}-\Ex{\ell_{k}}^2\\
		&=& \tfrac{1}{n-1}\tr\lp \Ex{L^2}\rp-(4m^2p^2)\lp\tfrac{1}{n-1}\rp^2 \\
		&=&r\sum_{i=1}^n\Ex{(L^2)_{ii}}-(4m^2p^2)r^2,
\end{eqnarray}
where the second equation follows
from the value of~$\mu$ and the fact
that the eigenvalues of~$L^2$ are the
squares of the eigenvalues of~$L$. 

We then expand to find
\begin{align}
\sigma^2 &= r\sum_{i=1}^n \left(\Ex{\lp\sum_{j=1, j\neq i}^n X_{ij}\rp^2}+\sum_{j=1,j\neq i}^n \Ex{X_{ij}^2}\right) -(4m^2p^2)r^2\\
		&= r\sum_{i=1}^n\lp pd_i +d_i(d_i-1)p^2+pd_i\rp-(4m^2p^2)r^2\label{eq:liisquared}, 
\end{align}
which follows from the independence of~$X_{ij}$ and~$X_{ik}$
(because edges are independent).
Re-arranging terms, we find
\begin{align}
\sigma^2 &= r\lp2p\sum_{i=1}^n d_i+p^2\sum_{i=1}^n d_i^2 \!-\! p^2\sum_{i=1}^n d_i\rp- 4m^2p^2r^2 \\
		&= r\lp 4mp+p^2\sum_{i=1}^n d_i^2 - 2p^2m\rp- 4m^2p^2r^2 \\
		&= \! \frac{p}{(n-1)^2}\lp  \!\! \lp n-1\rp \! \lp \! 2m(2-p) \!+\! p\sum_{i=1}^n  d_i^2\rp \!\!- 4m^2p\rp,
\end{align}
which follows from the Handshake Lemma. 
\hfill $\blacksquare$

With~$m, n \in \N$, 
$\underline{d}=(d_1,d_2,\dots d_n)$, and $0\le p \le 1$, we define
\begin{equation}
S(m,n,p,\underline{d}) = \left(2mp(n-1)(2-p) + p^2(n-1)\sum_{i=1}^n  d_i^2 -4m^2p^2\right)^{1/2}. 
\end{equation}
For the remainder of the paper, $\underline{d}$ will denote the degree sequence of $G$ (in any order).
A direct simplification of $\sigma^2$ above shows that it is equal to $\frac{1}{(n-1)^2}S(m,n,p,\underline{d})^2$.

\begin{lem}\label{lem:l1Nlowerbound}
Let~$G$ be connected and let~$\hat{G} \in \mathcal{G}(G, p)$.
Then, for~$\ell_{1:N}$ the smallest randomly sampled
eigenvalue of~$\hat{G}$, we have
\begin{equation}
\Ex{\ell_{1:N}} \geq \max\left\{0, \, \frac{2mp-S(m,n,p,\underline{d})\sqrt{N-1}}{n-1} \right\}. 
\end{equation}
\end{lem}

\emph{Proof:} Note that because any Laplacian is positive semi-definite, $\Ex{\ell_{1:N}}\geq 0$.
Then 
we apply Lemma \ref{lem:hale6} with $k=1$ to the values found in Lemma \ref{lem:musigma} to find that
\begin{align}
\Ex{\ell_{1:N}} &\geq \frac{2mp}{n-1} - \sigma \sqrt{\frac{N-1}{1}} \cdot \\
			&= \frac{1}{n-1}\lp 2mp-S(m,n,p,\underline{d})\sqrt{N-1} \rp.
\end{align}
\hfill $\blacksquare$

Next, we leverage this bound on~$\ell_{1:N}$ to bound
the expectation of~$\lambda_2$. 

\begin{lem}\label{lem:lambda2lowerbound}
Let~$G$ be connected and let~$\hat{G} \in \mathcal{G}(G, p)$,
where~$\hat{G}$ has Laplacian~$L$. For~$\lambda_2 := \lambda_2(L)$,
 we have 
\begin{equation}
\Ex{\lambda_2} \geq \max\left\{0, \left(1 \!\!-\!\! \left(\frac{n-2}{n-1}\right)^{N-1}\right)^{-1}\left[\frac{2mp}{n\!-\!1}\!\left(\!1 \!-\! \left(\!\frac{n-2}{n-1}\right)^{N-1}\right) \!-\! \frac{1}{n\!-\!1}S(m,n,p,\underline{d})\sqrt{N\!-\!1}\right]\!\right\}.
\end{equation}
\end{lem}
\emph{Proof:}
We put Lemma \ref{lem:l1Nupperbound} and Lemma \ref{lem:l1Nlowerbound} together to find that
\begin{equation}
\Ex{\lambda_2}\lp1-\lp\tfrac{n-2}{n-1}\rp^{N-1}\rp+\frac{2mp}{n-2}\lp \tfrac{n-2}{n-1}\rp^N \geq \frac{1}{n-1}\lp 2mp-S(m,n,p,\underline{d})\sqrt{N-1} \rp
\end{equation}
and solve this for $\Ex{\lambda_2}$.
\hfill $\blacksquare$

We will also require a bound on~$\lambda_2(L^2)$,
which we present next. 


\begin{lem}\label{lem:lambda2L2upperbound}
Let~$G$ be connected and let~$\hat{G} \in \mathcal{G}(G, p)$
have Laplacian~$L$. Then
\begin{equation}
\Ex{\lambda_2(L^2)}	\leq \frac{1}{n-1}\lp 4mp - 2mp^2+p^2\sum_{i=1}^n d_i^2\rp. 
\end{equation}
\end{lem}
\emph{Proof:} By definition, we
have~$\lambda_2(L^2) \leq \ell_k^2$ for an arbitrary~$k$.
Then
\begin{align}
\Ex{\lambda_2(L^2)}	&\leq \frac{1}{n-1}\tr{\Ex{L^2}}\\
					&= \frac{1}{n-1}\sum_{i=1}^n \Ex{\lp L^2\rp_{ii}}\\
					&= \frac{1}{n-1}\sum_{i=1}^n pd_i+d_i(d_i-1)p^2+pd_i,
\end{align}
which follows from Equation~\eqref{eq:liisquared}. 
Combining like terms, we find that					
\begin{align}
\Ex{\lambda_2(L^2)} &\leq \frac{1}{n-1}\lp\lp 2p\sum_{i=1}^n d_i \rp + \lp p^2\sum_{i=1}^n d_i^2\rp -\lp p^2\sum_{i=1}^n d_i \rp\rp\\
					&= \frac{1}{n-1}\lp 4mp - 2mp^2+p^2\sum_{i=1}^n d_i^2\rp,
\end{align}
as desired. 
\hfill $\blacksquare$

With these basic relationships involving~$\lambda_2$ now
established, the next section presents the solution
to Problem~\ref{prob:prob}.

\section{Probabilistic Connectivity of Random Graphs and Unions of Random Graphs} \label{sec:algconn}
This section translates the bounds on $\ell_{1:N}$ derived in Section~\ref{sec:firstorder}
for single random graphs into bounds on $\lambda_2$ for 
random graphs and unions of random
graphs. We then present our solution to 
Problem~\ref{prob:prob}.
In this section, we use the notation~$R(N, n)=1-((n-2)/(n-1))^{N-1}$.

\subsection{Solution to Problem \ref{prob:prob}}
Toward solving Problem~\ref{prob:prob}, we state
the Paley-Zygmund inequality. 

\begin{lem}[Paley-Zygmund inequality, \cite{paley32}]\label{lem:paleyzygmund}
Let $Z$ be a non-negative random variable with $\Var{Z} < \infty$ and let $\theta\in [0,1]$. Then
\begin{equation}
\Prob{Z>\theta\Ex{Z}} \geq \lp1-\theta\rp^2 \frac{\Ex{Z}^2}{\Ex{Z^2}}.
\end{equation}
\hfill $\blacksquare$
\end{lem}

We now present the main results of the paper,
namely a solution to Problem \ref{prob:prob}. 

\begin{thm} \label{thm:prob1} \emph{(Solution to Problem~\ref{prob:prob})} 
Let a connected graph~$G$ on~$n$ nodes be given, along
with~$p \in (0, 1)$. Let~$\hat{G} \in \mathcal{G}(G, p)$
and let~$L$ denote its Laplacian with second-smallest
eigenvalue~$\lambda_2$. 
Then
\begin{equation} \label{eq:thm1main}
\Prob{\lambda_2 > 0} \geq \max_{N \in \N} \frac{\lp \max\left\{\!0,\, 2mp R(N,n) -S(m,n,p,\underline{d})\sqrt{N - 1} \right\} \rp^2}{(n-1)R(N,n)^2 \lp 4mp - 2mp^2+p^2\sum_{i=1}^n d_i^2 \rp}.
\end{equation}
\end{thm}

\emph{Proof:}
We apply Lemma \ref{lem:paleyzygmund} with $\theta=0$ using the estimates for $\Ex{\lambda_2}$ from Lemma \ref{lem:lambda2lowerbound} and $\Ex{\lambda_2(L^2)}$ from Lemma \ref{lem:lambda2L2upperbound}. Because 
Equation~\eqref{eq:thm1main} holds for all~$N \in \N$,
it holds in particular for the maximum over~$N$. 
\hfill $\blacksquare$

Maximizing over $N$ is finitely terminating, as shown in the following remark.
\begin{rem}
Observe that for a given $\mathcal{G}(G,p)$, 
the parameter $N\in \N$ must satisfy
\begin{equation} 
1 \leq N \leq \frac{2mp(n-1)(2-p) + p^2(n-1)\sum_{i=1}^{n}d_i^2}{2mp(n-1)(2-p) + p^2(n-1)\sum_{i=1}^{n}d_i^2-4m^2p^2}
\end{equation}
to attain a positive value in the numerator of
Equation~\eqref{eq:thm1main}. Thus the maximum occurs
over this range of~$N$, and, because this range is finite,
searching for the maximum is finitely terminating. 
\hfill $\lozenge$
\end{rem}

Let~$\hat{p}(T) := 1 - (1 - p)^T$. We can then apply Theorem~\ref{thm:prob1} to graph unions in terms of the number of time steps required to attain connectivity. 

\begin{thm} \label{cor:prob1}
Given $\epsilon >0$, 
the graph~$\hat{G} \in \unionset_T(G,\hat{p}(T))$ is connected with probability $1-\epsilon$ after $T$ time steps if $T\geq T^*$, where $T^*$ is given by
\begin{equation}
T^* = \min\Bigg\{T \in \N : 
\max_{N \in \N} \!
\frac{\big[\max\{0,2m\hat{p}(T)R(N,n) - S(m,n,\hat{p}(T),\underline{d})\sqrt{N\!-\!1}\}\big]^2}
{(n\!-\!1)R(N,n)^2\big(4m\hat{p}(T) \!-\! 2m\hat{p}(T)^2 \!+\! \hat{p}(T)^2\sum_{i=1}^{n}d_i^2)} \geq 1 - \epsilon
\Bigg\}.
\end{equation}
\end{thm}
\emph{Proof:} Using Lemma~\ref{lem:unioneq}, this follows
from Theorem~\ref{thm:prob1} by replacing~$p$
with~$\hat{p}(T)$. \hfill $\blacksquare$

\begin{remark} \label{rem:T}
Examining Theorem~\ref{thm:prob1} and Theorem~\ref{cor:prob1},
we see that the probability of connectedness is monotonically
increasing in both~$p$ and~$T$. Thus, if agents in a 
network wish to 
increase the probability of connectedness of a union, they can
either increase~$p$, e.g., by spending more energy on
communications, or wait for~$T$ to increase as time elapses.
To quantify this relationship, 
recall from Lemma~\ref{lem:unioneq} that, in a union of~$T$
random graphs each with edge probability~$p$,
a single edge is present with 
probability~$\hat{p}(T) = 1 - (1 - p)^T$.
Thus, regardless of the value of~$p \in (0, 1)$, there is
some period of time~$T$ over which the probability
of an edge appearing can be made arbitrarily close to~$1$.
Consequently, varying~$T$ can drive the probability
of connectivity arbitrarily close to~$1$ as well,
regardless of the probability of connectivity of
any single graph in the union. \hfill $\Diamond$
\end{remark}

In the next section, we consider several specific
families of graphs and compute explicit lower
bounds for the solution to Problem~\ref{prob:prob} for these
families as~$n$ and~$p$ range across several orders
of magnitude. 

%

\section{Specific Graph Results} \label{sec:specific}
This section illustrates our main results for two
classes of underlying graphs:
complete graphs and complete graphs minus a cycle.

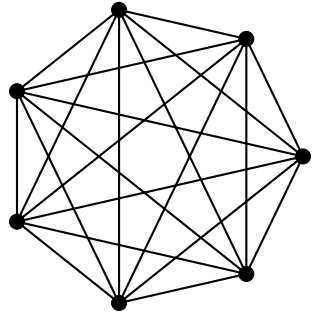
\begin{figure}
\centering
\begin{tikzpicture}[scale=0.5]
\foreach \a in {1,2,3,4,5,6,7}{
	\coordinate (\a) at (\a*360/7:4cm);
}
\draw[thick] (1) -- (2) -- (3) -- (4) -- (5) -- (6) -- (7) -- (1);
\draw[thick] (1) -- (3) -- (5) -- (7) -- (2) -- (4) -- (6) -- (1);
\draw[thick] (1) -- (4) -- (7) -- (3) -- (6) -- (2) -- (5) -- (1);
\foreach \b in {1,2,3,4,5,6,7}{
	\node[draw,circle,inner sep=2pt,fill] at (\b) {};
}	
\end{tikzpicture}
\caption{Complete graph on $7$ vertices, $K_7$}
\end{figure}

\subsection{Complete Graph}
We consider random graphs~$\hat{G} \in \mathcal{G}(K_n, p)$,
where~$K_n$ is the complete graph on~$n$ nodes and~$p \in (0, 1)$.
We note that this is equivalent to considering
conventional Erd\H{o}s-R\'{e}nyi random graphs and these results
therefore may be of independent interest. 

\begin{thm}
Let~$K_n$ denote the complete graph on~$n$ vertices
and let~$\hat{G} \in \mathcal{G}(K_n, p)$ have Laplacian~$L$
with second-smallest eigenvalue~$\lambda_2$. 
Then
\begin{equation}
\Prob{\lambda_2 > 0} \geq \\ \max_{N \in \N} \frac{\max\!\!\left\{\!0, \!\sqrt{n(n\!-\!1)p}R(\!N,n) \!-\!\! \sqrt{2(n\!-\!1)(1\!-\!p)(N\!-\!1)} \!\right\}^2}{(n-1)R(N,n)^2  \lp 2 - 2p+np \rp}.
\end{equation}
\end{thm}

\emph{Proof:}
We apply Theorem \ref{thm:prob1} with~$d_i=n-1$ for all~$i$
and~$m=\frac{n(n-1)}{2}$. 
Simplfying and factoring completes the proof.
\hfill $\blacksquare$

Similar to Theorem~\ref{cor:prob1}, 
we have the following corollary for unions
of graphs in~$\mathcal{G}(K_n, p)$. 
\begin{thm}
Let~$\hat{p} := \hat{p}(T)$. 
Given $\epsilon >0$, a graph~$\hat{U} \in \unionset_T(K_n,\hat{p}(T))$ is connected with probability $1-\epsilon$ after $T$ time steps if $T\geq T^*$, where $T^*$ is given by
\begin{equation}
T^* = \min \Bigg\{ T \in \N :  
 \max_{N \in \N} \frac{\max\!\!\left\{\!0, \!\sqrt{n(n\!-\!1)\hat{p}}R(\!N,n) \!-\!\! \sqrt{2(n\!-\!1)(1\!-\!\hat{p})(N\!-\!1)} \!\right\}^2}{(n-1)R(N,n)^2  \lp 2 - 2\hat{p}+n\hat{p} \rp} 
 \geq 1 - \epsilon \Bigg\}.
\end{equation}
\hfill $\blacksquare$
\end{thm}

\begin{figure}[t]
\centering
\includegraphics[width=3.3in]{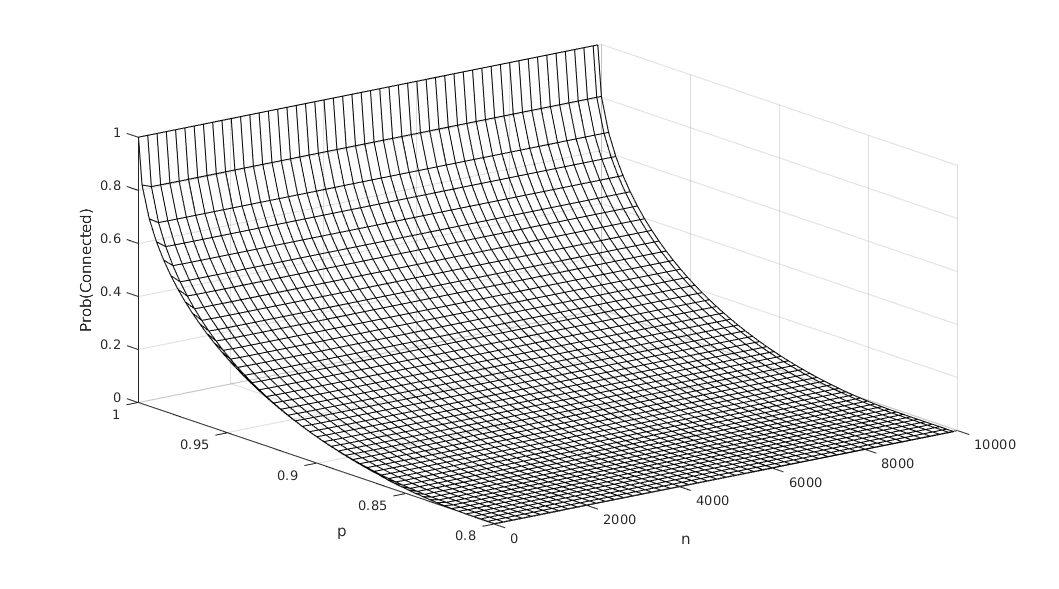}
\caption{Probability of connectedness of $\mathcal{G}(K_n,p)$
as a function of~$n$ and~$p$.
As noted in Remark~\ref{rem:T}, 
the desired probability of connectedness can be attained
for edge probabilities equal
to the depicted values of~$p$, or by waiting for~$T$
to be large enough that the value
of~$\hat{p}(T)  = 1 - (1 - p)^T$ is equal to the desired
depicted value of~$p$. 
}
\label{fig:Kn}
\end{figure}

We next examine numerical results. 
Figure~\ref{fig:Kn} shows the probability of connectedness
for a single random graph for~$1\leq n \leq 10,000$ 
and~$0.8\leq p \leq 1$. We see that increasing~$n$ provides
a modest increase in the probability of connectivity,
while increasing~$p$ makes connectivity substantially
more likely. In addition, as noted in Remark~\ref{rem:T},
accumulating more graphs in a union also increases
the probability of connectivity. In particular, even
if~$p$ is small, a large~$T$ can make~$\hat{p}(T)$ large,
thereby driving the probability of connectedness up
according to the increases in probabilities seen in
Figure~\ref{fig:Kn}.

\subsection{Complete Graph Minus a Cycle}
To illustrate our results beyond conventional Erd\H{o}s-R\'{e}nyi
graphs, we now consider random graphs in which
some edges never appear.
Formally, 
define $K_n\setminus C_n$ to be the complete graph on $n$ vertices less the edges of a cycle on $n$ vertices.  
Then we consider 
random graphs~$\hat{G} \in \mathcal{G}(K_n \setminus C_n, p)$.

\begin{figure}[ht]
\centering
\includegraphics[width=3.3in]{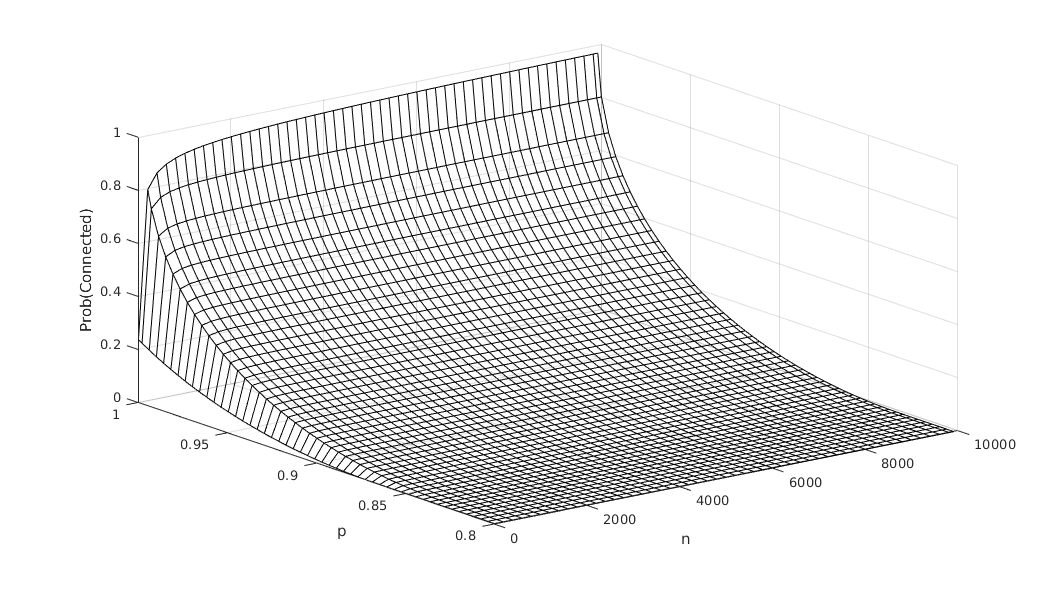}
\caption{Probability of connectedness of $\mathcal{G}(K_n\setminus C_n,p)$.
Here, as in Figure~\ref{fig:Kn}, making~$p$ larger for a single graph
(or making~$T$ larger to drive up~$\hat{p}(T)$ for a union of graphs)
increases the probability of connectivity. 
}
\label{fig:KnminusCn}
\end{figure}

Figure~\ref{fig:KnminusCn} shows the probability that 
a single graph~$\hat{G} \in \mathcal{G}(K_n\setminus C_n, p)$ is connected for $1\leq n \leq 10,000$ and $0.8\leq p \leq 1$. 
We find that, for a fixed pair~$(n, p)$, the probability
of connectivity is less than it is for the same~$(n,p)$
in Figure~\ref{fig:Kn}. This is intuitive, as 
a graph in~$\mathcal{G}(K_n\setminus C_n, p)$ has fewer
ways to attain connectivity than a graph
in~$\mathcal{G}(K_n, p)$. 

\begin{figure}[ht]
\centering
\begin{tikzpicture}[scale=0.5]
\foreach \a in {1,2,3,4,5,6,7}{
	\coordinate (\a) at (\a*360/7:4cm);
}
\draw[thick] (1) -- (3) -- (5) -- (7) -- (2) -- (4) -- (6) -- (1);
\draw[thick] (1) -- (4) -- (7) -- (3) -- (6) -- (2) -- (5) -- (1);
\foreach \b in {1,2,3,4,5,6,7}{
	\node[draw,circle,inner sep=2pt,fill] at (\b) {};
}	
\end{tikzpicture}
\caption{Complete graph on $7$ vertices with one cycle removed, $K_7 \setminus C_7$.}
\end{figure}
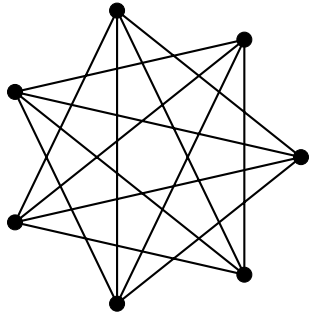

\section{Conclusion}\label{sec:conclusion}
We presented results bounding the probability of connectivity
for random graphs and unions thereof. These results can
be applied to multi-agent applications in both control
and optimization, where some convergence rates explicitly
depend upon the time needed to attain a connected union
graph. Future work includes extensions to time-varying
probabilities and heterogeneous edge probabilities.

\bibliographystyle{plain}{}
\bibliography{sources}

\section{Acknowledgements}
The authors would like to thank the following funding sources for enabling this research: the cooperative agreement FA8650-19-2-2429 with USRA (Unviersities Space Research Association); the STRESS (Science, Technology and Research for Exploiting Sensor Systems) contract (FA8650-18-2-1645); Air Force Office of Scientific Research grant FA9550-19-1-0169; US Department of Defense's Science, Mathematics and Research for Transformation (SMART) Scholarship for Service Program; and the Air Force Research Lab Sensors Directorate.

\end{document}